\documentclass[12pt,reqno]{amsart}
\usepackage[utf8]{inputenc}
\usepackage{fancyhdr}
\usepackage{amssymb}
\usepackage{amsthm}
\usepackage{amsmath}
\usepackage{mathtools}
\usepackage{amsxtra}
\usepackage{mathrsfs}
\usepackage[all]{xy}
\usepackage{colonequals}
\usepackage{enumitem}
\usepackage{fullpage}
\usepackage{breqn}
\usepackage{caption}
\usepackage{placeins}

\numberwithin{equation}{section}
\theoremstyle{plain}
\newtheorem{thm}{Theorem}[section]
\newtheorem{lem}[thm]{Lemma} 
\newtheorem{prop}[thm]{Proposition}
\newtheorem{cor}[thm]{Corollary}

\theoremstyle{definition}
\newtheorem*{ack*}{Acknowledgment}
\theoremstyle{remark}

\newtheorem{rmk}{Remark}[section]

\DeclareMathOperator{\li}{li}

\newcommand{\C}{\mathbb C}

\newcommand{\R}{\mathbb R}

\input cyracc.def

\usepackage{amsrefs}
\usepackage[hidelinks]{hyperref}
\makeatletter
\@namedef{subjclassname@2020}{\textup{2020} Mathematics Subject Classification}
\makeatother

\begin{document}
\title{Numerically Explicit Estimates for the Distribution of Rough Numbers}
\author{Kai (Steve) Fan}
\address{Department of Mathematics, Dartmouth College, Hanover, NH 03755, USA}
\email{steve.fan.gr@dartmouth.edu}
\setlength{\footskip}{8mm}
\begin{abstract}
For $x\ge y>1$ and $u\colonequals \log x/\log y$, let $\Phi(x,y)$ denote the number of positive integers up to $x$ free of prime divisors less than or equal to $y$. In 1950 de Bruijn \cite{Br} studied the approximation of $\Phi(x,y)$ by the quantity
\[\mu_y(u)e^{\gamma}x\log y\prod_{p\leq y}\left(1-\frac{1}{p}\right),\]
where $\gamma=0.5772156...$ is Euler's constant and 
\[\mu_y(u)\colonequals\int_{1}^{u}y^{t-u}\omega(t)\,dt.\]
He showed that the asymptotic formula
\[\Phi(x,y)=\mu_y(u)e^{\gamma}x\log y\prod_{p\leq y}\left(1-\frac{1}{p}\right)+O\left(\frac{xR(y)}{\log y}\right)\]
holds uniformly for all $x\ge y\ge2$, where $R(y)$ is a positive decreasing function related to
the error estimates in the Prime Number Theorem. In this paper we obtain numerically explicit versions of de Bruijn's result.
\end{abstract}

\maketitle
\section{Introduction}\label{S:Intro}
Let $x\ge y>1$ be positive real numbers. Throughout the paper, we shall always write $u\colonequals \log x/\log y$, and the letters $p$ and $q$ will always denote primes. We say that a positive integer $n$ is \emph{$y$-rough} if all the prime divisors of $n$ are greater than $y$. Let $\Phi(x,y)$ denote the number of $y$-rough numbers up to $x$. Explicitly, we have
\[\Phi(x,y)=\sum_{\substack{n\leq x\\P^-(n)>y}}1,\]
where $P^-(n)$ denotes the least prime divisor of $n$, with the convention that $P^-(1)=\infty$. When $1\le u\le 2$, or equivalently when $\sqrt{x}\le y\le x$, we simply have $\Phi(x,y)=\pi(x)-\pi(y)+1$, where $\pi(\cdot)$ is the prime-counting function. The function $\Phi(x,y)$ is closely related to the sieve of Eratosthenes, one of the most ancient algorithms for finding primes, and it has been extensively studied by mathematicians. A simple application of the inclusion-exclusion principle enables us to write
\begin{equation}\label{Equ:IE}
\Phi(x,y)=\sum_{d\mid P(y)}\mu(d)\left\lfloor\frac{x}{d}\right\rfloor,
\end{equation}
where $\lfloor a\rfloor $ is the integer part of $a$ for any $a\in\R$, $\mu$ is the M\"{o}bius function, and $P(y)$ denotes the product of primes up to $y$. If $y$ is relatively small in comparison with $x$, say $y=x^{o(1)}$, the above formula can be used to obtain $\Phi(x,y)\sim e^{-\gamma}x/\log y$ as $y\to\infty$, where $\gamma=0.5772156...$ is Euler's constant. However, it turns out that this nice asymptotic formula does not hold uniformly, as already exemplified by the base case $1\le u\le 2$.


In 1937, Buchstab \cite{Bu} showed that for any fixed $u>1$, one has $\Phi(x,y)\sim \omega(u)x/\log y$ as $x\to\infty$, where $\omega(u)$ is defined to be the unique continuous solution to the delay differential equation $(u\omega(u))'=\omega(u-1)$ for $u\ge2$, subject to the initial value condition $\omega(u)=1/u$ for $u\in[1,2]$. Comparing this result with the asymptotic formula obtained from (\ref{Equ:IE}), one would expect that $\omega(u)\to e^{-\gamma}$ as $u\to\infty$. Indeed, it can be shown \cite[Corollary III.6.5]{T} that $\omega(u)=e^{-\gamma}+O(u^{-u/2})$ for $u\ge1$. Moreover, it is known that $\omega(u)$ oscillates above and below $e^{-\gamma}$ infinitely often. It is convenient to extend the definition of $\omega(u)$ by setting $\omega(u)=0$ for all $u<1$, so that $\omega(u)$ satisfies the same delay differential equation on $\R\setminus\{1,2\}$. In the sequel, we shall write $\omega'(1)$ and $\omega'(2)$ for the right derivatives of $\omega(u)$ at $u=1$ and $u=2$, respectively. With this convention, we have $(u\omega(u))'=\omega(u-1)$ for all $u\in\R$.
\par Buchstab's asymptotic formula can be proved easily based on the following identity \cite[Theorem III.6.3]{T} named after him:
\begin{equation}\label{Equ:BuchstabEqu}
\Phi(x,y)=\Phi(x,z)+\sum_{y<p\le z}\sum_{v\ge1}\Phi(x/p^{v},p)
\end{equation}
for any $z\in[y,x]$. The Buchstab function $\omega(u)$ then appears naturally in the iteration process, starting with $\Phi(x,y)\sim x/(u\log y)$ in the range $1< u\le 2$. Since $1/2\le\omega(u)\le1$ for $u\in[1,\infty)$, Buchstab's asymptotic formula suggests that the relation $\Phi(x,y) \asymp x/\log y$ holds uniformly for $x\ge y>1$. Thus, it is of interest to seek numerically explicit estimates for $\Phi(x,y)$ that are applicable in wide ranges. Confirming a conjecture of Ford, the author \cite{F} showed that $\Phi(x,y)<x/\log y$ holds uniformly for $x\ge y>1$, which is essentially best possible when $x^{1-\epsilon}\le y\le \epsilon x$, where $\epsilon\in(0,1)$ is fixed. On the other hand, the values of $\omega(u)$ indicate that improvements should be expected in the narrower range $2\le y\le\sqrt{x}$. In recent work jointly with Pomerance \cite{FP}, the author proved that $\Phi(x,y)<0.6x/\log y$ holds uniformly for $3\le y\le\sqrt{x}$. This inequality provides a fairly good upper bound for $\Phi(x,y)$, especially considering that the absolute maximum of $\omega(u)$ over $[2,\infty)$ is given by $M_0=0.5671432...$, attained at the unique critical point $u=2.7632228...$ of the function $(\log(u-1)+1)u^{-1}$ on $[2,3]$. With a bit more effort, one can show, using the Buchstab identity (\ref{Equ:BuchstabEqu}), that
\begin{equation}\label{Equ:implicitDelta(x,y)}
\Phi(x,y)=\frac{x}{\log y}\left(\omega(u)+O\left(\frac{1}{\log y}\right)\right)
\end{equation}
uniformly for $2\le y\le\sqrt{x}$ (see \cite[Theorem III.6.4]{T}). In Section \ref{S:0.4Inequ}, we shall derive a numerically explicit lower bound of this type that suits our needs. Our method can also be modified with ease to obtain a numerically explicit upper bound of the same type. 
\par In \cite{Br} de Bruijn provided a more precise approximation for $\Phi(x,y)$ than $\omega(u)x/\log y$. Let us fix some $y_0\ge2$ for the moment. Suppose that there exist a positive constant $C_0(y_0)$ and a positive decreasing function $R(z)$ defined on $[y_0,\infty)$, such that $R(z)\gg\, z^{-1}$, that $R(z)\downarrow 0$ as $z\to\infty$ and that for all $z\ge y_0$ we have
\begin{equation}\label{Equ:pi(x)}
|\pi(z)-\li(z)|\le\frac{z}{\log z}R(z)
\end{equation}
and
\begin{equation}\label{Equ:intpi(x)}
\int_{z}^{\infty}\frac{|\pi(t)-\li(t)|}{t^2}\,dt\le C_0(y_0)R(z),
\end{equation}
where $\li(z)$ is the logarithmic integral defined by
\[\li(z)\colonequals\int_{0}^{z}\frac{dt}{\log t}.\]
The classical version of the Prime Number Theorem allows us to take $R(z)=\exp(-c\sqrt{\log z})$ for some suitable constant $c>0$. Using the zero-free region of Korobov and Vinogradov for the Riemann zeta-function, we obtain $R(z)=\exp(-c'(\log z)^{3/5}(\log\log z)^{-1/5})$ for some absolute constant $c'>0$. If the Riemann Hypothesis holds, then one can take $R(z)=c''z^{-1/2}\log^2z$, where $c''>0$ is an absolute constant.
\par To state de Bruijn's result, we define 
\[\mu_y(u)\colonequals\int_{1}^{u}y^{t-u}\omega(t)\,dt.\]
It is easy to see that $0\le\mu_y(u)\log y\le1-y^{1-u}$ and that for every fixed $u\ge1$, we have $\mu_y(u)\log y\to\omega(u)$ as $y\to\infty$. Precise expansions for $\mu_y(u)$ in terms of the powers of $\log y$ can be found in \cite[Theorem III.6.18]{T}. When $1\le u\le 2$, the change of variable $t=\log v/\log y$ shows that 
\[\mu_y(u)x=\int_{1}^{u}t^{-1}y^t\,dt=\int_{y}^{x}\frac{dv}{\log v}=\li(x)-\li(y).\]
Since $\Phi(x,y)=\pi(x)-\pi(y)+1$ when $1\le u\le 2$, (\ref{Equ:pi(x)}) clearly implies that
\[\Phi(x,y)=\mu_y(u)x+(\pi(x)-\li(x))-(\pi(y)-\li(y))+1=\mu_y(u)x+O\left(\frac{xR(y)}{\log y}\right).\]
It can be shown using (\ref{Equ:pi(x)}) and (\ref{Equ:intpi(x)}) that
\[\prod_{p\leq y}\left(1-\frac{1}{p}\right)=\frac{e^{-\gamma}}{\log y}\left(1+O(R(y))\right).\]
Thus we have, equivalently,
\begin{equation}\label{Equ:de Bruijn}
	\Phi(x,y)=\mu_y(u)e^{\gamma}x\log y\prod_{p\leq y}\left(1-\frac{1}{p}\right)+O\left(\frac{xR(y)}{\log y}\right).
\end{equation}
Essentially, de Bruijn \cite{Br} showed that this formula holds uniformly for $x\ge y\ge y_0$. In Section \ref{S:de Bruijn} we shall derive an explicit version of (\ref{Equ:de Bruijn}), which will be applied in Section \ref{S: Pf of Thm2&Cor3} to obtain numerically explicit estimates with suitable $y_0$ and $R(y)$. Our main results are summarized in the following theorem.
\begin{thm}\label{thm:explicit debruijn}
Uniformly for $x\ge y\ge 2$, we have
\[\left|\Phi(x,y)-\mu_y(u)e^{\gamma}x\log y\prod_{p\leq y}\left(1-\frac{1}{p}\right)\right|< 4.403611\frac{x}{(\log y)^{3/4}}\exp\left(-\sqrt{\frac{\log y}{6.315}}\right).\]
Conditionally on the Riemann Hypothesis, we have
\[\left|\Phi(x,y)-\mu_y(u)e^{\gamma}x\log y\prod_{p\leq y}\left(1-\frac{1}{p}\right)\right|<0.449774\frac{x\log y}{\sqrt{y}}\]
uniformly for $x\ge y\ge 11$. 
\end{thm}
The following consequence of Theorem \ref{thm:explicit debruijn} is sometimes more convenient to use.
\begin{cor}\label{cor::explicit debruijn}
Uniformly for $x\ge y\ge 2$, we have
\[|\Phi(x,y)-\mu_y(u)x|< 4.434084\frac{x}{(\log y)^{3/4}}\exp\left(-\sqrt{\frac{\log y}{6.315}}\right).\]
Conditionally on the Riemann Hypothesis, we have
\[|\Phi(x,y)-\mu_y(u)x|< 0.460680\frac{x\log y}{\sqrt{y}}\]
uniformly for $x\ge y\ge 11$. 
\end{cor}

\vspace*{2mm}
\section{Lower Bounds for $\Phi(x,y)$}\label{S:0.4Inequ}
Before moving on to the derivation of Theorem \ref{thm:explicit debruijn}, we prove a clean lower bound for $\Phi(x,y)$ which is applicable in a wide range. This lower bound, which is interesting in itself, will be used in the proof of Theorem \ref{thm:explicit debruijn} and Corollary \ref{cor::explicit debruijn} in Section \ref{S: Pf of Thm2&Cor3}. We start by proving the following result, which provides a numerically explicit lower bound for the implicit constant in the error term in (\ref{Equ:implicitDelta(x,y)}). As we already mentioned, our method can easily be adapted to yield a numerically explicit upper bound as well, though it will not be needed in the present paper.
\begin{prop}\label{Prop:LBDelta(x,y)}
Define $\Delta(x,y)$ by
\[\Phi(x,y)=\frac{x}{\log y}\left(\omega(u)+\frac{\Delta(x,y)}{\log y}\right)\]
for $2\le y\le \sqrt{x}$. Let $y_0=602$. For every positive integer $k\ge3$, we define
\[\Delta_k^-=\Delta_k^-(y_0)\colonequals\inf\left\{\min(\Delta(x,y),0)\colon y\ge y_0\emph{~and~}2\le u<k\right\}.\]
Then $\Delta_3^->-0.563528$, $\Delta_4^->-0.887161$, and $\Delta_k^->-0.955421$ for all $k\ge5$.
\end{prop}
\begin{proof}
Let $y_1\colonequals2{,}278{,383}$. Suppose first that $y\ge y_1$ and set
\[G(v)\colonequals\sum_{x^{1/v}<p\le \sqrt{x}}\frac{1}{p}\]
for $2\le v\le u$. By \cite[Theorem 5.6]{D}\footnote{In \cite{B} it is claimed that the proof of  \cite[Theorem 4.2]{D} is incorrect due to the application of an incorrect zero density estimate of Rama\'re \cite[Theorem 1.1]{Rama}. In a footnote on p. 2299 of the same paper, the authors state that the bounds asserted in \cite{D} are likely affected for this reason. However, since they also give a correct proof of \cite[Theorem 4.2]{D} (see \cite[Corollary 11.2]{B}), one verifies easily that the proof of \cite[Theorem 5.6]{D}, which relies only on \cite[Theorem 4.2]{D}, partial summation, and numerical computation, remains valid.}, we have
\begin{equation}\label{PropEqu1}
\left|G(v)-\log\frac{v}{2}\right|\le\frac{c_1}{\log^2 y}
\end{equation}
for all $y\ge y_1$, where $c_1=0.4/\log y_1$. We shall also make use of the following inequality \cite[Corollary 5.2]{D}\footnote{For the same reason mentioned above, it is reasonable to suspect that the bounds given in \cite[Corollary 5.2]{D} are also affected. However, one can verify these bounds without much difficulty. Indeed, (5.2) of \cite[Corollary 5.2]{D} is superseded by \cite[Corollary 1]{RS}, while (5.3) and (5.4) of \cite[Corollary 5.2]{D} follow from \cite[Lemmas 3.2--3.4]{BD} and direct calculations.}:
\begin{equation}\label{PropEqu2}
\frac{z}{\log z}\left(1+\frac{c_3}{\log z}\right)\le\pi(z)\le \frac{z}{\log z}\left(1+\frac{c_2}{\log z}\right),
\end{equation} 
where $c_2=1+2.53816/\log y_1$ and $c_3=1+2/\log y_1$. We start with the range $2\le u\le3$. In this range, we have
\begin{align*}
\Phi(x,y)&=\#\{n\le x\colon P^-(n)>y\text{~and~}\Omega(n)\le2\}\\
&=\pi(x)-\pi(y)+1+ \sum_{y<p\leq \sqrt{x}}\sum_{p\le q\le x/p}1\\
&=\pi(x)-\pi(y)+1+ \sum_{y<p\leq \sqrt{x}}(\pi(x/p)-\pi(p)+1),
\end{align*}
where $\Omega(n)$ denotes the total number of prime factors of $n$, with multiplicity counted. 
Since
\[\sum_{y<p\le \sqrt{x}}(\pi(p)-1)=\sum_{\pi(y)<j\le\pi(\sqrt{x})}(j-1)=\frac{\pi(\sqrt{x})(\pi(\sqrt{x})-1)}{2}-\frac{\pi(y)(\pi(y)-1)}{2},\]
we see that
\[\pi(x)-\pi(y)+1-\sum_{y<p\le \sqrt{x}}(\pi(p)-1)>\pi(x)-\frac{\pi(\sqrt{x})^2}{2}+\frac{\pi(\sqrt{x})}{2}.\]
It follows from (\ref{PropEqu2}) that 
\begin{equation}\label{PropEqu3}
\Phi(x,y)>\frac{x}{\log x}\left(1+\frac{c_3}{\log x}\right)-\frac{x}{2\log^2\sqrt{x}}\left(1+\frac{c_2}{\log\sqrt{x}}\right)^2+\frac{\sqrt{x}}{2\log\sqrt{x}}+\sum_{y<p\leq \sqrt{x}}\pi(x/p).
\end{equation}
To handle the sum in (\ref{PropEqu3}), we appeal to (\ref{PropEqu2}) again to arrive at
\[\sum_{y<p\leq \sqrt{x}}\pi(x/p)\ge\sum_{y<p\le\sqrt{x}}\left(\frac{x}{p\log(x/p)}+\frac{c_3x}{p\log^2(x/p)}\right).\]
By partial summation we see that
\[\sum_{y<p\le \sqrt{x}}\frac{1}{p\log(x/p)}=\frac{1}{\log x}\int_{2^-}^{u}\frac{v}{v-1}\,dG(v)=\frac{1}{\log y}\left(\frac{G(u)}{u-1}+\frac{1}{u}\int_{2^-}^{u}\frac{G(v)}{(v-1)^2}\,dv\right)\]
From (\ref{PropEqu1}) it follows that
\[\frac{G(u)}{u-1}\ge\frac{1}{u-1}\left(\log\frac{u}{2}-\frac{c_1}{\log^2y}\right),\]
and
\begin{align*}
\int_{2^-}^{u}\frac{G(v)}{(v-1)^2}\,dv&\ge\int_{2}^{u}\frac{1}{(v-1)^2}\left(\log\frac{v}{2}-\frac{c_1}{\log^2y}\right)\,dv\\
&=-\frac{1}{u-1}\log\frac{u}{2}+\int_{2}^{u}\frac{1}{v(v-1)}\,dv-\frac{c_1}{\log^2y}\left(1-\frac{1}{u-1}\right)\\
&=-\frac{u}{u-1}\log\frac{u}{2}+\log(u-1)-\frac{c_1}{\log^2y}\left(1-\frac{1}{u-1}\right).
\end{align*}
Hence
\begin{align}\label{PropEqu4}
\sum_{y<p\le \sqrt{x}}\frac{x}{p\log(x/p)}&\ge \frac{x}{\log y}\left(\frac{\log(u-1)}{u}-\frac{2c_1}{u\log^2y}\right) \nonumber\\
&=\frac{x}{\log y}\left(\omega(u)-\frac{2c_1}{u\log^2y}\right)-\frac{x}{\log x}.
\end{align}
Similarly, we have
\[\sum_{y<p\le \sqrt{x}}\frac{1}{p\log^2(x/p)}=\frac{1}{\log^2 x}\int_{2^-}^{u}\left(\frac{v}{v-1}\right)^2\,dG(v)=\frac{1}{\log^2 x}\left(\frac{G(u)u^2}{(u-1)^2}+2\int_{2^-}^{u}\frac{vG(v)}{(v-1)^3}\,dv\right).\]
By (\ref{PropEqu1}) we have
\[\frac{G(u)u^2}{(u-1)^2}\ge\frac{u^2}{(u-1)^2}\left(\log\frac{u}{2}-\frac{c_1}{\log^2y}\right),\]
and
\[\int_{2^-}^{u}\frac{vG(v)}{(v-1)^3}\,dv\ge\int_{2}^{u}\frac{v}{(v-1)^3}\left(\log\frac{v}{2}-\frac{c_1}{\log^2y}\right)\,dv.\]
Since
\begin{align*}
\int_{2}^{u}\frac{v}{(v-1)^3}\log\frac{v}{2}\,dv&=-\left(\frac{1}{u-1}+\frac{1}{2(u-1)^2}\right)\log\frac{u}{2}+\int_{2}^{u}\left(\frac{1}{v-1}+\frac{1}{2(v-1)^2}\right)\frac{dv}{v}\\
&=-\frac{2u-1}{2(u-1)^2}\log\frac{u}{2}+\frac{1}{2}\int_{2}^{u}\left(\frac{1}{(v-1)^2}+\frac{1}{v(v-1)}\right)\,dv\\
&=-\frac{u^2}{2(u-1)^2}\log\frac{u}{2}+\frac{1}{2}\left(\log(u-1)+1-\frac{1}{u-1}\right)
\end{align*}
and 
\[\int_{2}^{u}\frac{v}{(v-1)^3}\,dv=-\frac{2u-1}{2(u-1)^2}+\frac{3}{2},\]
we have
\begin{equation}\label{PropEqu5}
\sum_{y<p\le \sqrt{x}}\frac{x}{p\log^2(x/p)}\ge\frac{x}{\log^2 x}\left(\log(u-1)+\frac{u-2}{u-1}-\frac{4c_1}{\log^2y}\right).
\end{equation}
Inserting (\ref{PropEqu4}) and (\ref{PropEqu5}) into (\ref{PropEqu3}) yields
\[\Delta(x,y)\ge g(u)-\frac{2c_1}{u\log y}+\frac{\log y}{uy^{3/2}}-\frac{1}{u^2}\left(2-c_3+\frac{4c_1c_3}{\log^2 y}+\frac{8c_2}{u\log y}+\frac{8c_2^2}{u^2\log^2 y}\right),\]
where
\[g(u)\colonequals \frac{c_3}{u^2}\left(\log(u-1)+\frac{u-2}{u-1}\right).\]
Using Mathematica we find that $\Delta_3^->-0.301223$ when $y\ge y_1$.
\par Now we proceed to  bound $\Delta_k^-$ for $k\ge4$ recursively when $y\ge y_1$. Let $k\ge3$ be arbitrary. It is easily seen that the following variant of Buchstab's identity (\ref{Equ:BuchstabEqu}) holds for any $z\in[y,x]$:
\begin{equation}\label{Equ:VBuchstab}
\Phi(x,y)=\Phi(x,z)+\sum_{y<p\le z}\Phi(x/p,p^-),
\end{equation}
where $p^-<p$ is any real number sufficiently close to $p$. For $3\le k\le u<k+1$ and $y\ge y_1$, we obtain by taking $z=x^{1/3}$ that
\begin{equation}\label{PropEqu6}
\Phi(x,y)=\Phi\left(x,x^{1/3}\right)+\sum_{y<p\le x^{1/3}}\Phi(x/p,p^-).
\end{equation}
We have already shown that
\begin{equation}\label{PropEqu7}
\Phi\left(x,x^{1/3}\right)\ge\frac{x}{\log x^{1/3}}\left(\omega\left(\frac{\log x}{\log x^{1/3}}\right)+\frac{\Delta_3^-}{\log x^{1/3}}\right)=\frac{3x}{\log y}\left(\frac{\omega(3)}{u}+\frac{3\Delta_3^-}{u^2\log y}\right).
\end{equation}
Note that $2<\log(x/p)/\log(p^-)<k$. Thus, we have 
\[\Phi(x/p,p^-)\ge\frac{x}{p\log(p^-)}\left(\omega\left(\frac{\log(x/p)}{\log(p^-)}\right)+\frac{\Delta_k^-}{\log(p^-)}\right).\]
Since $\omega(u)$ is continuous on $[1,\infty)$, it follows from (\ref{PropEqu6}) and (\ref{PropEqu7}) that
\begin{equation}\label{PropEqu8}
\Phi(x,y)\ge\frac{3x}{\log y}\left(\frac{\omega(3)}{u}+\frac{3\Delta_3^-}{u^2\log y}\right)+\sum_{y<p\le x^{1/3}}\frac{x}{p\log p}\left(\omega\left(\frac{\log x}{\log p}-1\right)+\frac{\Delta_k^-}{\log p}\right).
\end{equation}
By partial summation we see that
\begin{equation*}
\sum_{y<p\le x^{1/3}}\frac{1}{p\log^2 p}<\int_{y}^{\infty}\frac{1}{t\log^2t}\,d\pi(t)=-\frac{\pi(y)}{y\log^2y}+\int_{y}^{\infty}\frac{\log t+2}{t^2\log^3t}\pi(t)\,dt,
\end{equation*}
which, by (\ref{PropEqu2}), is
\begin{align*}
&<-\frac{1}{\log^3y}\left(1+\frac{c_3}{\log y}\right)+\int_{y}^{\infty}\frac{\log t+2}{t\log^4t}\left(1+\frac{c_2}{\log t}\right)\,dt \\
&=-\frac{1}{\log^3y}\left(1+\frac{c_3}{\log y}\right)+\frac{1}{2\log^2y}+\frac{c_2+2}{3\log^3y}+\frac{c_2}{2\log^4y}\\
&=\frac{1}{\log^2y}\left(\frac{1}{2}+\left(\frac{c_2}{3}-1\right)\frac{1}{\log y}+\left(\frac{c_2}{2}-c_3\right)\frac{1}{\log^2y}\right)\\
&<\frac{1}{\log^2y}\left(\frac{1}{2}+\left(\frac{c_2}{3}-1\right)\frac{1}{\log y}\right).
\end{align*}
Hence
\begin{equation}\label{PropEqu9}
\sum_{y<p\le x^{1/3}}\frac{\Delta_k^-x}{p\log^2 p}\ge\frac{\Delta_k^-x}{\log^2y}\left(\frac{1}{2}+\left(\frac{c_2}{3}-1\right)\frac{1}{\log y}\right).
\end{equation}
On the other hand, we have
\begin{align*}
\sum_{y<p\le x^{1/3}}\frac{1}{p\log p}\omega\left(\frac{\log x}{\log p}-1\right)&=\frac{1}{\log x}\int_{3^-}^{u}v\omega(v-1)\,dG(v)\\
&=\frac{1}{\log x}\left(\int_{3}^{u}\omega(v-1)\,dv+\int_{3^-}^{u}v\omega(v-1)\,d\left(G(v)-\log\frac{v}{2}\right)\right).
\end{align*}
Observe that 
\[\int_{3}^{u}\omega(v-1)\,dv=u\omega(u)-3\omega(3)\]
and that
\begin{align*}
\int_{3^-}^{u}v\omega(v-1)\,d\left(G(v)-\log\frac{v}{2}\right)&=u\omega(u-1)\left(G(v)-\log\frac{v}{2}\right)-3\omega(2)\left(G(3)-\log\frac{3}{2}\right)\\
	&\hspace*{3mm}-\int_{3^-}^{u}\left(G(v)-\log\frac{v}{2}\right)\,d(v\omega(v-1)).
\end{align*}
By \cite[(6.23), p. 562]{T} and \cite[Theorems III.5.7 \& III.6.6]{T}, we have, for all $v\ge3$, that
\begin{align*}
	\frac{d}{dv}(v\omega(v-1))&=\omega(v-2)+\omega'(v-1)\ge\frac{1}{2}-\rho(v-1)\ge\frac{1}{2}-\rho(2)=\log2-\frac{1}{2},\\
	\frac{d}{dv}(v\omega(v-1))&\le 1+\rho(v-1)\le1+\rho(2)=2-\log2,
\end{align*}
where $\rho$ is the Dickman-de Bruijn function defined to be the unique continuous solution to the delay differential equation $t\rho'(t)+\rho(t-1)=0$ for $t\ge1$, subject to the initial value condition $\rho(t)=1$ for $0\le t\leq 1$. Moreover, we have
\[\lim\limits_{v\to3^-}\frac{d}{dv}(v\omega(v-1))=\lim\limits_{v\to3^-}\left(\omega(v-2)+\omega'(v-1)\right)=-\frac{1}{4}.\]
It follows by (\ref{PropEqu1}) that
\[\int_{3^-}^{u}\left(G(v)-\log\frac{v}{2}\right)\,d(v\omega(v-1))\le \frac{c_1}{\log^2y}\left(u\omega(u-1)-3\omega(2)\right).\]
Thus we have
\[\int_{3^-}^{u}v\omega(v-1)\,d\left(G(v)-\log\frac{v}{2}\right)\ge -\frac{2c_1u\omega(u-1)}{\log^2y}\ge -\frac{2c_1M_0u}{\log^2y},\]
where $M_0=0.5671432...$. Hence we have shown that
\begin{equation}\label{PropEqu10}
	\sum_{y<p\le x^{1/3}}\frac{x}{p\log p}\omega\left(\frac{\log x}{\log p}-1\right)\ge\frac{x}{\log y}\left(\omega(u)-\frac{3\omega(3)}{u}- \frac{2c_1M_0}{\log^2 y}\right).
\end{equation}
Combining (\ref{PropEqu8}), (\ref{PropEqu9}) and (\ref{PropEqu10}), we deduce that
\[\Delta(x,y)\ge\frac{9\Delta_3^-}{u^2}+\frac{\Delta_k^-}{2}-\frac{1}{\log y}\left(2c_1M_0-\left(\frac{c_2}{3}-1\right)\Delta_k^-\right)\]
for $k\le u<k+1$. Therefore, $\Delta_{k+1}^-\ge \min(\Delta_k^-,a_k^-)$ for all $k\ge3$, where
\[a_k^-\colonequals\frac{9\Delta_3^-}{k^2}+\frac{\Delta_k^-}{2}-\frac{1}{\log y_1}\cdot\max\left(2c_1M_0-\left(\frac{c_2}{3}-1\right)\Delta_k^-,0\right).\]
Consequently, we have $\Delta_4^->-0.451835$ and $\Delta_k^->-0.480075$ for all $k\ge5$.
\par Suppose now that $602\le y\le y_1$. By \cite[Theorem 20]{RS} we can replace (\ref{PropEqu1}) with
\[\left|G(v)-\log\frac{v}{2}\right|\le\frac{d_1}{\sqrt{y}\log y},\]
where $d_1=2$. Moreover, (\ref{PropEqu2}) remains true if we replace $c_2$ and $c_3$ by $d_2=1.2762$ and $d_3=1$, respectively, according to \cite[Corollary 5.2]{D}. With these changes, we run the same argument used to handle the case $y\ge y_1$ and get
\[\Delta(x,y)>g(u)-\frac{2d_1}{u\sqrt{y}}+\frac{\log y}{uy^{3/2}}-\frac{1}{u^2}\left(2-d_3+\frac{4d_1d_3}{\sqrt{y}\log y}+\frac{8d_2}{u\log y}+\frac{8d_2^2}{u^2\log^2 y}\right).\]
when $2\le u\le 3$ and 
\[\Delta(x,y)\ge\frac{9\Delta_3^-}{u^2}+\frac{\Delta_k^-}{2}-\frac{1}{\log y}\left(\frac{2d_1M_0\log y}{\sqrt{y}}-\left(\frac{d_2}{3}-1\right)\Delta_k^-\right)\]
when $3\le k\le u<k+1$, so that we can take
\[a_k^-=\frac{9\Delta_3^-}{k^2}+\frac{\Delta_k^-}{2}-\frac{1}{\log y_0}\cdot\max\left(\frac{2d_1M_0\log y_0}{\sqrt{y_0}}-\left(\frac{d_2}{3}-1\right)\Delta_k^-,0\right).\]
As a consequence, we have $\Delta_3^->-0.563528$, $\Delta_4^->-0.887161$ and $\Delta_k^->-0.955421$ for all $k\ge5$. This completes the proof of the proposition.
\end{proof}

The next result provides a numerical lower bound for $\omega(u)$ on $[3,\infty)$.
\begin{lem}\label{Lem:minomega(u)foru>=3}
We have $\omega(u)>0.549307$ for all $u\ge3$.
\end{lem}
\begin{proof}
Consider first the case $u\in[3,4]$. Since $(t\omega(t))'=\omega(t-1)$ for $t\ge2$ and $\omega(t)=(\log(t-1)+1)/t$ for $t\in[2,3]$, we have
\[\omega(u)=\frac{1}{u}\left(\log 2+1+\int_{3}^{u}\frac{\log(t-2)+1}{t-1}\,dt\right)\]
for $u\in[3,4]$. Note that $u\omega'(u)=\omega(u-1)-\omega(u)=S(u)/u$, where
\[S(u)\colonequals\frac{u(\log(u-2)+1)}{u-1}-\log2-1-\int_{3}^{u}\frac{\log(t-2)+1}{t-1}\,dt.\]
Since
\begin{align*}
S'(u)&=\frac{1}{u-1}\left(\log(u-2)+1+\frac{u}{u-2}-\frac{u(\log(u-2)+1)}{u-1}-(\log(u-2)+1)\right)\\
&=\frac{u(1-(u-2)\log(u-2))}{(u-2)(u-1)^2},
\end{align*}
we know that $S(u)$ is strictly increasing on $[3,u_1]$ and strictly decreasing on $[u_1,4]$, where $u_1=3.7632228...$ is the unique solution to the equation $(u-2)\log(u-2)=1$. But $S(3)=1/2-\log2<0$ and
\[S(4)=\frac{\log2+1}{3}-\int_{3}^{4}\frac{\log(t-2)+1}{t-1}\,dt>0.\]
It follows that $S(u)$ has a unique zero $u_2\in[3,4]$. The numerical value of $u_2$ is given by $u_2=3.4697488...$, according to Mathematica. Hence $S(u)<0$ for $u\in[3,u_2)$ and $S(u)>0$ for $u\in(u_2,4]$. The same is true for $\omega'(u)$, which implies that $\omega(u)$ is strictly decreasing on $[3,u_2]$ and strictly increasing on $[u_2,4]$. Thus, $\omega(u)\ge\omega(u_2)=0.5608228...$ for $u\in[3,4]$. 
\par Consider now the case $u\in[4,\infty)$. It is known \cite{JR} that $\omega(t)$ satisfies
\[|\omega(t)-e^{-\gamma}|\leq \frac{\rho(t-1)}{t}\]
for all $t\geq1$. Since $\rho(t)$ is strictly decreasing on $[4,\infty)$, we have $\omega(u)\ge e^{-\gamma}-\rho(3)/4$ for all $u\ge4$. To find the value of $\rho(3)$, we use $t\rho'(t)+\rho(t-1)=0$ for $t\ge1$ and $\rho(t)=1-\log t$ for $t\in[1,2]$ to obtain
\[\rho(u)=1-\log 2-\int_{2}^{u}\frac{1-\log(t-1)}{t}\,dt\]
for $u\in[2,3]$. It follows that
\[\omega(u)\ge e^{-\gamma}-\frac{1}{4}\left(1-\log 2-\int_{2}^{3}\frac{1-\log(t-1)}{t}\,dt\right)=0.5493073...\]
for all $u\ge4$. We have therefore shown that $\omega(u)>0.549307$ for all $u\ge3$.
\end{proof}

We are now ready to prove the following clean lower bound for $\Phi(x,y)$ that we alluded to.
\begin{thm}\label{thm:0.4LB}
We have $\Phi(x,y)>0.4x/\log y$ uniformly for all $7\le y\le x^{2/3}$.
\end{thm}
\begin{proof}
In the range $\max(7,x^{2/5})\le y\le x^{2/3}$, we have trivially $\Phi(x,y)\ge \pi(x)-\pi(y)+1$. By \cite[Corollary 5.2]{D} we have
\begin{align*}
\pi(x)-\pi(y)&\ge\frac{x}{\log x}\left(1+\frac{1}{\log x}\right)-\frac{y}{\log y}\left(1+\frac{1.2762}{\log y}\right)\\
&=\left(\frac{1}{u}\left(1+\frac{1}{\log x}\right)-\frac{y}{x}\left(1+\frac{1.2762u}{\log x}\right)\right)\frac{x}{\log y}\\
&>\left(\frac{2}{5}\left(1+\frac{1}{\log x}\right)-\frac{1}{x^{1/3}}\left(1+\frac{3.1905}{\log x}\right)\right)\frac{x}{\log y}>0.4\frac{x}{\log y}
\end{align*}
whenever $x\ge41{,}217$. Furthermore, we have verified $\Phi(x,y)>0.4x/\log y$ for $\max(7,x^{2/5})\le y\le x^{2/3}$ with $x\le41{,}217$ using Mathematica. Hence, $\Phi(x,y)>0.4x/\log y$ holds in the range $\max(7,x^{2/5})\le y\le x^{2/3}$.
\par Consider now the case $\max(x^{1/3},7)\le y\le x^{2/5}$. Following the proof of Proposition \ref{Prop:LBDelta(x,y)}, we have
\begin{align}\label{Inequ:PhiM}
\Phi(x,y)&=\pi(x)-\pi(y)+1+ \sum_{y<p\leq \sqrt{x}}(\pi(x/p)-\pi(p)+1)\nonumber\\
&=\pi(x)-M(x,y)+\sum_{y<p\le x^{1/2}}\pi(x/p),
\end{align}
where
\[M(x,y)\colonequals\frac{1}{2}\pi\left(\sqrt{x}\right)^2-\frac{1}{2}\pi\left(\sqrt{x}\right)-\frac{1}{2}\pi(y)^2+\frac{3}{2}\pi(y)-1.\]
To handle the sum in (\ref*{Inequ:PhiM}), we appeal to Theorem 5 and its corollary from \cite{RS} to arrive at
\[G(v)-\log\frac{v}{2}>-\frac{1}{2\log^2\sqrt{x}}-\frac{1}{\log^2 y}\ge-\frac{33}{25\log^2y}\]
in the range $\max(x^{1/3},7)\le y\le x^{2/5}$. By \cite[Corollary 1]{RS} we have
\[\sum_{y<p\le x^{1/2}}\pi(x/p)>x\sum_{y<p\le x^{1/2}}\frac{1}{p\log(x/p)}=\frac{x}{\log x}\int_{2^-}^{u}\frac{v}{v-1}\,dG(v),\]
provided that $x\ge289$. The right-hand side of the above can be estimated in the same way as in the proof of Proposition \ref{Prop:LBDelta(x,y)}, so we obtain
\[\sum_{y<p\le \sqrt{x}}\pi(x/p)> \frac{x}{\log y}\left(\omega(u)-\frac{66}{25u\log^2y}\right)-\frac{x}{\log x}.\]
On the other hand, we see by \cite[Corollary 5.2]{D} and \cite[Corollary 2]{RS} that
\[\pi(x)-M(x,y)>\pi(x)-\frac{1}{2}\pi\left(\sqrt{x}\right)^2\ge\frac{x}{\log x}\left(1+\frac{1}{\log x}\right)-\frac{25x}{8\log^2x}=\frac{x}{\log x}-\frac{17x}{8\log^2x}.\]
for $x\ge114^2$. Collecting the estimates above and using the inequality $\omega(u)\ge\omega(5/2)=2(\ln(3/2)+1)/5$ for $u\in[5/2,3]$, we find that
\[\Phi(x,y)>\frac{\omega(5/2)x}{\log y}-\frac{17x}{8\log^2x}-\frac{66x}{25u\log^3y}\ge \frac{\omega(5/2)x}{\log y}-\frac{17x}{50\log^2y}-\frac{132x}{125\log^3y}>0.4\frac{x}{\log y}\]
for all $\max(46,x^{1/3})\le y\le x^{2/5}$. For $x^{1/3}\le y\le x^{2/5}$ with $7\le y\le 46$, we have verified the inequality $\Phi(x,y)>0.4x/\log y$ directly through numerical computation.
\par Next, we consider the range $7\le y<x^{1/3}$. By Proposition \ref{Prop:LBDelta(x,y)} and Lemma \ref{Lem:minomega(u)foru>=3} we have
\[\Phi(x,y)>\frac{x}{\log y}\left(0.549307-\frac{0.955421}{\log y}\right)>0.4\frac{x}{\log y},\]
provided that $y\ge 602$. To deal with the range $7\le y\le \min(x^{1/3},602)$, we follow the inclusion-exclusion technique used in \cite[Section 3]{FP}. For any integer $n\ge1$, let $\nu(n)$ denote the number of distinct prime factors of $n$. We start by ``pre-sieving" with the primes 2, 3, and 5: for any $x\ge1$ the number of integers
$n\le x$ with $\gcd(n,30)=1$ is $(4/15)x+r_x$, where $|r_x|\le14/15$. Let $P_5(y)$ be the product of the primes in $(5,y]$.  Then we have by the Bonferroni inequalities that
\[\Phi(x,y)\ge\sum_{\substack{d\mid P_5(y)\\\nu(d)\le3}}\mu(d)\left(\frac{4}{15}\cdot\frac{x}{d}+r_{x/d}\right)\ge a(y)x-b(y),\]
where
\begin{align*}
a(y)&\colonequals\frac{4}{15}\sum_{\substack{d\mid P_5(y)\\\nu(d)\le3}}\frac{\mu(d)}{d}=\frac{4}{15}\sum_{j=0}^3(-1)^j\sum_{\substack{d\mid P_5(y)\\\nu(d)=j}}\frac{1}{d},\\
b(y)&\colonequals\frac{14}{15}\sum_{j=0}^{3}\binom{\pi(y)-3}{j}.
\end{align*}
By Newton's identities, the inner sum in the definition of $a(y)$ can be represented in terms of the power sums of $1/p$ over all primes $5<p\le y$. Thus, we have $\Phi(x,y)>0.4x/\log y$ whenever $a(y)>0.4/\log y$ and $x> b(y)/(a(y)-0.4/\log y)$. Using Mathematica, we find that the inequality $\Phi(x,y)>0.4x/\log y$ holds for $7\le y\le 602$ and $x\ge 13{,}160{,}748$. Finally, we have verified the inequality $\Phi(x,y)>0.4x/\log y$ directly for $7\le y\le x^{1/3}$ with $x\le13{,}160{,}748$ by numerical calculations, completing the proof of our theorem.
\end{proof}

\begin{rmk}
Note that for $y\in[5,7)$ we have
\[\Phi(x,y)\ge\frac{4}{15}x-\frac{14}{15}>0.4\frac{x}{\log 5}\ge0.4\frac{x}{\log y},\]
provided that $x\ge 52$. Combined with Theorem \ref{thm:0.4LB} and numerical examination of the case $11\le x\le 52$, this implies that the inequality $\Phi(x,y)>0.4x/\log y$ holds in the slightly larger range $5\le y\le x^{2/3}$ if one assumes $x\ge 41$.
\end{rmk}

\vspace*{2mm}
\section{An Explicit Version of de Bruijn's Estimate}\label{S:de Bruijn}
To prove Theorem \ref{thm:explicit debruijn}, we shall first develop an explicit version of (\ref{Equ:de Bruijn}) with a general $R(y)$, following \cite{Br}, where $R(y)$ is a positive decreasing function satisfying the same conditions described in the introduction. Suppose that $y_0\ge3$. For each $z\ge 2$, put
\[Q(z)\colonequals\prod_{p\leq z}\left(1-\frac{1}{p}\right).\]
We start by estimating $Q(y)$ for $y\ge y_0$. Using a Stieltjes integral, we may write
\begin{equation}\label{Equ1:logQ(z)/Q(y)}
\log\frac{Q(z)}{Q(y)}=\int_{y}^{z}\log\left(1-t^{-1}\right)\,d\li(y)+\int_{y}^{z}\log\left(1-t^{-1}\right)\,d(\pi(y)-\li(t)),
\end{equation}
where $z\ge y\ge y_0$. The first integral on the right-hand side of the above is equal to
\begin{align*}
\int_{y}^{z}\log\left(1-t^{-1}\right)\frac{dt}{\log t}=-\log\frac{\log z}{\log y}+\int_{y}^{z}\left(t^{-1}+\log\left(1-t^{-1}\right)\right)\frac{dt}{\log t}.
\end{align*}
Since 
\[ -\frac{1}{2t(t-1)}<t^{-1}+\log\left(1-t^{-1}\right)<0\]
for all $t\ge y_0$, we have 
\[-\frac{1}{2}\int_{y}^{\infty}\frac{dt}{t(t-1)\log t}<\int_{y}^{z}\left(t^{-1}+\log\left(1-t^{-1}\right)\right)\frac{dt}{\log t}<0.\]
But a change of variable shows that
\[\int_{y}^{\infty}\frac{dt}{t(t-1)\log t}=\int_{1}^{\infty}\frac{dt}{t(y^t-1)}\le\frac{1}{y-1}\int_{1}^{\infty}\frac{dt}{t^2}=\frac{1}{y-1},\]
where we have used the inequality $y^t-1\ge (y-1)t$ for $t\ge 1$ and $y\ge y_0$. It follows that
\begin{equation}\label{Equ:1stint}
-\frac{1}{2(y-1)}\le\int_{y}^{z}\log\left(1-t^{-1}\right)\,d\li(y)+\log\frac{\log z}{\log y}<0.
\end{equation}
Now we estimate the second integral on the right-hand side of (\ref{Equ1:logQ(z)/Q(y)}). By (\ref{Equ:pi(x)}) and partial integration we have
\begin{align*}
\left|\int_{y}^{z}\log\left(1-t^{-1}\right)\,d(\pi(y)-\li(t))\right|&\le\log\left(1-y^{-1}\right)^{-1}\frac{y}{\log y}R(y)+\log\left(1-z^{-1}\right)^{-1}\frac{z}{\log z}R(z)\\
&\hspace*{2mm}+\int_{y}^{z}\frac{|\pi(t)-\li(t)|}{t(t-1)}\,dt.
\end{align*}
Using (\ref{Equ:intpi(x)}) we see that 
\[\int_{y}^{z}\frac{|\pi(t)-\li(t)|}{t(t-1)}\,dt\le \frac{C_0(y_0)y_0}{y_0-1}R(y).\]
It is clear that the function 
\[\log\left(1-t^{-1}\right)^{-1}\frac{t}{\log t}=\frac{1}{\log t}\sum_{n=0}^{\infty}\frac{t^{-n}}{n+1}\]
is strictly decreasing for $t\in(1,\infty)$. Since $R(t)$ is decreasing on $[y_0,\infty)$, we find that
\[\left|\int_{y}^{z}\log\left(1-t^{-1}\right)\,d(\pi(y)-\li(t))\right|\le \left(2\log\left(1-y_0^{-1}\right)^{-1}\frac{y_0}{\log y_0}+\frac{C_0(y_0)y_0}{y_0-1}\right)R(y).\]
Combining this inequality with (\ref{Equ1:logQ(z)/Q(y)}) and (\ref{Equ:1stint}) yields
\begin{equation}\label{Equ2:logQ(z)/Q(y)}
-C_2(y_0)R(y)\le\log\frac{Q(z)}{Q(y)}+\log\frac{\log z}{\log y}\le C_1(y_0)R(y)
\end{equation}
for $z\ge y\ge y_0$, where 
\begin{align*}
C_1(y_0)&=2\log\left(1-y_0^{-1}\right)^{-1}\frac{y_0}{\log y_0}+\frac{C_0(y_0)y_0}{y_0-1},\\
C_2(y_0)&=C_1(y_0)+\sup_{t\ge y_0}\frac{1}{2(t-1)R(t)}.
\end{align*}
Exponentiating (\ref{Equ2:logQ(z)/Q(y)}) we obtain
\begin{equation}\label{Equ:Q(z)/Q(y)}
-C_4(y_0)R(y)\le\frac{Q(z)\log z}{Q(y)\log y}-1\le C_3(y_0)R(y)
\end{equation}
for $z\ge y\ge y_0$, where
\begin{align*}
C_3(y_0)&=\sup_{t\ge y_0}\frac{\exp(C_1(y_0)R(t))-1}{R(t)}=\frac{\exp(C_1(y_0)R(y_0))-1}{R(y_0)},\\
C_4(y_0)&=\sup_{t\ge y_0}\frac{1-\exp(-C_2(y_0)R(t))}{R(t)}=C_2(y_0).
\end{align*}
As a consequence, we have by letting $z\to\infty$ in (\ref{Equ:Q(z)/Q(y)}) and using the fact that $Q(z)\log z\to e^{-\gamma}$ as $z\to\infty$, that
\begin{equation}\label{Equ:1/Q(y)}
e^{\gamma}\log y(1-C_4(y_0)R(y))\le \frac{1}{Q(y)}\le e^{\gamma}\log y(1+C_3(y_0)R(y)).
\end{equation}
Similarly, we derive from (\ref{Equ2:logQ(z)/Q(y)}) that
\begin{equation}\label{Equ:Q(y)}
\frac{e^{-\gamma}}{\log y}(1-C_6(y_0)R(y))\le Q(y)\le \frac{e^{-\gamma}}{\log y}(1+C_5(y_0)R(y))
\end{equation}
for $y\ge y_0$, where
\begin{align*}
C_5(y_0)&=\sup_{t\ge y_0}\frac{\exp(C_2(y_0)R(t))-1}{R(t)}=\frac{\exp(C_2(y_0)R(y_0))-1}{R(y_0)},\\
C_6(y_0)&=\sup_{t\ge y_0}\frac{1-\exp(-C_1(y_0)R(t))}{R(t)}=C_1(y_0).
\end{align*}
\par For $x\ge y\ge 2$, we define 
\[\psi(x,y)\colonequals \frac{\Phi(x,y)}{xQ(y)}.\]
We then need to estimate $\eta(x,y)=\psi(x,y)-\lambda(x,y)$, where $\lambda(x,y)\colonequals e^{\gamma}\mu_y(u)\log y$. For $1\le u\le 2$ this can be done straightforward. Indeed, we have $\Phi(x,y)=\pi(x)-\pi(y)+1$ and $\omega(u)=1/u$ when $1\le u\le 2$, so that
\[\eta(x,y)=\frac{\pi(x)-\pi(y)+1}{xQ(y)}-e^{\gamma}\log y\int_{1}^{u}t^{-1}y^{t-u}\,dt.\]
Note that
\[\left|\pi(x)-\pi(y)-x\int_{1}^{u}t^{-1}y^{t-u}\,dt\right|=\left|\pi(x)-\pi(y)-\int_{y}^{x}\frac{dt}{\log t}\right|\le \left(\frac{x}{\log x}+\frac{y}{\log y}\right)R(y).\]
From (\ref{Equ:1/Q(y)}) it follows that $|\eta(x,y)|\le e^{\gamma}\alpha_y(u)R(y)$ for $y\ge y_0$ and $u\in[1,2]$, where
\[\alpha_y(u)\colonequals\frac{\log y}{y^uR(y)}+C_3(y_0)\left(\frac{\log y}{y^u}+\log y\int_{1}^{u}t^{-1}y^{t-u}\,dt\right)+(1+C_3(y_0)R(y))\left(\frac{1}{u}+y^{1-u}\right).\]
Integration by parts enables us to write
\begin{equation*}
\log y\int_{1}^{u}t^{-1}y^{t-u}\,dt=\frac{1}{u}-y^{1-u}+\int_{1}^{u}t^{-2}y^{t-u}\,dt
\end{equation*}
for $y\ge y_0$. Hence $|\eta(x,y)|\le e^{\gamma}\eta_1(y)R(y)$ for $y\ge y_0$ and $u\in[1,2]$, where
\begin{equation}\label{Equ:eta_1}
\eta_1(y)\colonequals\sup_{t\ge y}\frac{\log t}{tR(t)}+\max_{u\in[1,2]}\left(C_3(y_0)I_{y}(u)+(1+C_3(y_0)R(y))\left(\frac{1}{u}+y^{1-u}\right)\right)
\end{equation}
with 
\[I_{y}(u)\colonequals \frac{1}{u}+\int_{1}^{u}t^{-2}y^{t-u}\,dt.\]
We remark that $I_{y}(u)$ is strictly decreasing on $[1,2]$ and hence satisfies $I_{y}(u)<1$ for $u\in(1,2]$, since its derivative is
\[I_{y}'(u)=-\int_{1}^{u}t^{-2}y^{t-u}\log y\,dt<0.\]
Thus, (\ref{Equ:eta_1}) simplifies to 
\begin{equation}\label{Equ:eta_1'}
\eta_1(y)=\sup_{t\ge y}\frac{\log t}{tR(t)}+C_3(y_0)+2(1+C_3(y_0)R(y)).
\end{equation}
\par Suppose now that $y\ge y_0$ and $u\ge2$. From (\ref{Equ:VBuchstab}) it follows that
\begin{equation}\label{Equ:psi1}
\psi(x,y)=\psi(x,z)\frac{Q(z)}{Q(y)}+\sum_{y<p\le z}\psi(x/p,p^-)\cdot\frac{1}{p}\prod_{y<q<p}\left(1-\frac{1}{q}\right),
\end{equation}
where $z\ge y\ge y_0$. Put $h\colonequals \log z/\log y\ge1$ and
\begin{equation}\label{Equ:H_y(v)}
H_y(v)\colonequals\sum_{y<p\le y^v}\frac{1}{p}\prod_{y<q<p}\left(1-\frac{1}{q}\right)
\end{equation}
for $v\ge1$. Then we have $H_y(v)=1-Q(y^v)/Q(y)$. By partial summation, we see that (\ref{Equ:psi1}) becomes
\begin{equation}\label{Equ:psi2}
\psi(x,y)=\psi(y^u,y^h)(1-H_y(h))+\int_{1}^{h}\psi(y^{u-v},(y^v)^-)\,dH_y(v).
\end{equation}
By (\ref{Equ:Q(z)/Q(y)}) we have
\[|H_y(v)-1+v^{-1}|\le C_7(y_0)R(y),\]
where $C_7(y_0)\colonequals\max(C_3(y_0),C_4(y_0))$. Thus, one can think of $1-v^{-1}$ as a smooth approximation to $H_y(v)$. Since we also expect $\lambda(x,y)$ to be a smooth approximation to $\psi(x,y)$, in view of (\ref{Equ:psi2}) it is reasonable to expect 
\[E_1(h;y,u)\colonequals\lambda(y^u,y)-\lambda(y^u,y^h)h^{-1}-\int_{1}^{h}\lambda(y^{u-v},y^v)v^{-2}\,dv\]
to be small in size as a function of $y$. This can be easily verified when $1\le h\le u/2$. Following de Bruijn \cite{Br}, we have 
\begin{equation}\label{Equ:partialE_1}
\frac{\partial}{\partial h}E_1(h;y,u)=-h^{-1}\cdot\frac{\partial}{\partial h}\lambda(y^u,y^h)+h^{-2}\lambda(y^{u},y)-h^{-2}\lambda(y^{u-h},y^h).
\end{equation}
Since 
\[\frac{\lambda(y^u,y^h)}{e^{\gamma}\log y}=h\int_{1}^{u/h}y^{ht-u}\omega(t)\,dt,\]
we find
\[\frac{\partial}{\partial h}\left(\frac{\lambda(y^u,y^h)}{e^{\gamma}\log y}\right)=\int_{1}^{u/h}y^{ht-u}\omega(t)\,dt+h\left(\log y\int_{1}^{u/h}y^{ht-u}(t\omega(t))\,dt-uh^{-2}\omega(uh^{-1})\right).\]
Recall that $(t\omega(t))'=\omega(t-1)$ for $t\in\R$ with the obvious extension $\omega(t)=0$ for $t<1$. It follows that
\begin{align*}
\log y\int_{1}^{u/h}y^{ht-u}(t\omega(t))\,dt&=h^{-1}y^{ht-u}(t\omega(t))\bigg|_{1}^{u/h}-h^{-1}\int_{1}^{u/h}y^{ht-u}\omega(t-1)\,dt\\
&=uh^{-2}\omega(uh^{-1})-h^{-1}y^{h-u}-h^{-1}y^h\int_{1}^{u/h-1}y^{ht-u}\omega(t)\,dt\\
&=uh^{-2}\omega(uh^{-1})-h^{-1}y^{h-u}-\left(h^2e^\gamma\log y\right)^{-1}\lambda(y^{u-h},y^h).
\end{align*}
Hence we have
\begin{align*}
\frac{\partial}{\partial h}\lambda(y^u,y^h)&=e^{\gamma}\log y\left(\int_{1}^{u/h}y^{ht-u}\omega(t)\,dt-y^{h-u}\right)-h^{-1}\lambda(y^{u-h},y^h)\\
&=h^{-1}\lambda(y^u,y^h)-e^{\gamma}y^{h-u}\log y-h^{-1}\lambda(y^{u-h},y^h).
\end{align*}
Inserting this in (\ref{Equ:partialE_1}) yields
\[\frac{\partial}{\partial h}E_1(h;y,u)=h^{-1}e^{\gamma}y^{h-u}\log y.\]
Integrating both sides with respect to $h$ and using the initial value condition $E_1(1;y,u)=0$, we obtain 
\begin{equation}\label{Equ:E_1}
E_1(h;y,u)=e^{\gamma}\log y\int_{1}^{h}t^{-1}y^{t-u}\,dt<e^{\gamma}y^{h-u}.
\end{equation}
\par In what follows, we shall always suppose that $1\le h\le u/2$. Following de Bruijn \cite{Br}, we proceed to show that
\[E_3(h;y,u)\colonequals \lambda(y^u,y)-\lambda(y^u,y^h)(1-H(h))-\int_{1}^{h}\lambda(y^{u-v},y^v)\,dH(h)\]
is small in size as a function of $y$. This is intuitive, because
\[\lambda(y^u,y^h)h^{-1}-\int_{1}^{h}\lambda(y^{u-v},y^v)v^{-2}\,dv,\]
which is a good approximation to $\lambda(y^u,y)$ as we have already demonstrated, can be thought of as a smooth approximation to
\[\lambda(y^u,y^h)(1-H(h))-\int_{1}^{h}\lambda(y^{u-v},y^v)\,dH(h).\]
Moreover, we have by (\ref{Equ:psi2}) that
\begin{equation}\label{Equ:eta(x,y)1}
\eta(x,y)=\eta(y^u,y^h)(1-H_y(h))+\int_{1}^{h}\eta(y^{u-v},(y^v)^-)\,dH_y(v)-E_3(h;y,u),
\end{equation}
which will later be used to estimate $\eta(x,y)$. To estimate $E_3(h;y,u)$, let us write $E_3(h;y,u)=E_1(h;y,u)+E_2(h;y,u)$, where 
\[E_2(h;y,u)\colonequals-\int_{1}^{h}\lambda(y^{u-v},y^v)\,d\left(H(v)-1+v^{-1}\right)+(H(h)-1+h^{-1})\lambda(y^u,y^h).\]Then we expect $E_2(h;y,u)$ to be small in size as a function of $y$. Using (\ref{Equ:H_y(v)}) and the observation that $H(1)=0$, we have
\begin{equation}\label{Equ:E_2}
|E_2(h;y,u)|\le \left(\left|\lambda(y^u,y^h)-\lambda(y^{u-h},y^h)\right|+\int_{1}^{h}\left|\frac{\partial}{\partial v}\lambda(y^{u-v},y^v)\right|\,dv\right)C_7(y_0)R(y).
\end{equation}
Note that
\begin{align*}
\frac{\lambda(y^u,y^h)-\lambda(y^{u-h},y^h)}{he^{\gamma}\log y}&=\int_{1}^{u/h}y^{ht-u}\omega(t)\,dt-\int_{2}^{u/h}y^{ht-u}\omega(t-1)\,dt\\
&=\int_{1}^{2}y^{ht-u}\omega(t)\,dt+\int_{2}^{u/h}y^{ht-u}(\omega(t)-\omega(t-1))\,dt\\
&=\int_{1}^{2}t^{-1}y^{ht-u}\,dt-\int_{2}^{u/h}y^{ht-u}t\omega'(t)\,dt.
\end{align*}
By Theorems III.5.7 and III.6.6 in \cite{T} we have
\begin{equation}\label{Equ:omega'}
|\omega'(t)|\le\rho(t)\le\frac{1}{\Gamma(t+1)}
\end{equation}
for all $t\ge1$. It follows that
\begin{equation}\label{Equ:lamba-lambda}
\left|\lambda(y^u,y^h)-\lambda(y^{u-h},y^h)\right|\le he^{\gamma}\log y\left(\int_{1}^{2}t^{-1}y^{ht-u}\,dt+\int_{2}^{u/h}y^{ht-u}t\rho(t)\,dt\right).
\end{equation}
This inequality will later be used in conjunction with the formulas
\begin{equation}\label{Equ:int1}
h\log y\int_{1}^{2}t^{-1}y^{ht-u}\,dt=\frac{y^{2h-u}}{2}-y^{h-u}+\int_{1}^{2}t^{-2}y^{ht-u}\,dt
\end{equation}
and 
\begin{align}
h\log y\int_{2}^{u/h}y^{ht-u}t\rho(t)\,dt&=uh^{-1}\rho(uh^{-1})-2\rho(2)y^{2h-u}-\int_{2}^{u/h}y^{ht-u}(t\rho(t))'\,dt\nonumber\\
&\le uh^{-1}\rho(uh^{-1})-2\rho(2)y^{2h-u}+\int_{2}^{u/h}y^{ht-u}\rho(t-1)\,dt.\label{Equ:int2}
\end{align}
On the other hand, we have
\[\frac{\lambda(y^{u-v},y^v)}{e^{\gamma}\log y}=v\int_{2}^{u/v}y^{vt-u}\omega(t-1)\,dt,\]
which implies that
\[\frac{\partial}{\partial v}\left(\frac{\lambda(y^{u-v},y^v)}{e^{\gamma}\log y}\right)=\int_{2}^{u/v}y^{vt-u}(1+tv\log y)\omega(t-1)\,dt-uv^{-1}\omega(uv^{-1}-1).\]
By partial integration, the right side of the above is easily seen to be
\[-2y^{2v-u}-\int_{2}^{u/v}y^{vt-u}t\omega'(t-1)\,dt.\]
Hence, we arrive at
\[\int_{1}^{h}\left|\frac{\partial}{\partial v}\lambda(y^{u-v},y^v)\right|\,dv\le e^{\gamma}\log y\left(2\int_{1}^{h}y^{2v-u}\,dv+\int_{1}^{h}\int_{2}^{u/v}y^{vt-u}t|\omega'(t-1)|\,dtdv\right).\]
Furthermore, we have by Fubini's theorem that
\[\int_{1}^{h}\int_{2}^{u/v}y^{vt-u}t|\omega'(t-1)|\,dtdv=\int_{2}^{u/h}\int_{1}^{h}y^{vt-u}t|\omega'(t-1)|\,dv\,dt+\int_{u/h}^{u}\int_{1}^{u/t}y^{vt-u}t|\omega'(t-1)|\,dv\,dt,\]
the right side of which is easily seen to be
\[\frac{1}{\log y}\left(\int_{2}^{u/h}y^{ht-u}|\omega'(t-1)|\,dt+\int_{u/h}^{u}|\omega'(t-1)|\,dt-\int_{2}^{u}y^{t-u}|\omega'(t-1)|\,dt\right).\]
It follows that
\begin{equation}\label{Equ:partiallambda}
\int_{1}^{h}\left|\frac{\partial}{\partial v}\lambda(y^{u-v},y^v)\right|\,dv< e^{\gamma}\left(y^{2h-u}+\int_{2}^{u/h}y^{ht-u}|\omega'(t-1)|\,dt+\int_{u/h}^{u}|\omega'(t-1)|\,dt\right).
\end{equation}
This estimate together with (\ref*{Equ:lamba-lambda}) will lead us to a good estimate for $E_2(h;y,u)$. 
\par Now we derive estimates for $E_3(h;y,u)$ that suit our needs. Suppose that $k\le u< k+1$ and take $h=h_k=u/k$, where $k\ge2$ is a positive integer. We first consider the case $k=2$. In view of (\ref{Equ:int1}), we see that (\ref{Equ:lamba-lambda}) simplifies to
\[\left|\lambda\left(y^u,y^{h_2}\right)-\lambda\left(y^{u-h_2},y^{h_2}\right)\right|< e^{\gamma}\left(\frac{1}{2}+\int_{1}^{2}t^{-2}y_0^{t-2}\,dt\right)=e^{\gamma}I_{y_0}(2)\]
for $y\ge y_0$. By (\ref*{Equ:partiallambda}) we have
\[\int_{1}^{h_2}\left|\frac{\partial}{\partial v}\lambda(y^{u-v},y^v)\right|\,dv\le e^{\gamma}\left(1+\int_{2}^{3}|\omega'(t-1)|\,dt\right)=\frac{3e^{\gamma}}{2},\]
since $\omega'(t)=-1/t^2$ for $t\in[1,2)$. Combining these estimates with (\ref{Equ:E_1}) and (\ref{Equ:E_2}), we obtain $E_3(h_2;y,u)\le e^{\gamma}\xi_2(y_0)R(y)$ for $y\ge y_0$ and $2\le u<3$, where
\[\xi_2(y_0)\colonequals \max_{t\ge y_0}\frac{1}{tR(t)}+C_7(y_0)\left(I_{y_0}(2)+\frac{3}{2}\right).\]
Now we handle the case $k\ge3$. From (\ref{Equ:omega'})--(\ref*{Equ:int2}) it follows that
\begin{dmath*}
\left|\lambda\left(y^u,y^{h_k}\right)-\lambda\left(y^{u-h_k},y^{h_k}\right)\right|< e^{\gamma}\left(\frac{1}{\Gamma(k)}+\left(2\log 2-\frac{3}{2}\right)y^{2-k}+\int_{1}^{2}t^{-2}y^{t-k}\,dt+\int_{2}^{3}y^{t-k}(1-\log(t-1))\,dt+\int_{3}^{k}y^{t-k}\frac{dt}{\Gamma(t)}\right),
\end{dmath*}
where we have used the fact that $\rho(t)=1-\log t$ for $t\in[1,2]$. By (\ref{Equ:omega'}) and (\ref{Equ:partiallambda}) we have
\[\int_{1}^{h_k}\left|\frac{\partial}{\partial v}\lambda(y^{u-v},y^v)\right|\,dv\le e^{\gamma}\left(y^{2-k}+\int_{2}^{3}y^{t-k}\frac{dt}{(t-1)^2}+\int_{3}^{k}y^{t-k}\frac{dt}{\Gamma(t)}+\int_{k}^{k+1}\frac{dt}{\Gamma(t)}\right).\]
Together with (\ref{Equ:E_1}) and (\ref{Equ:E_2}), these inequalities imply that $E_3(h_k;y,u)\le e^{\gamma}\xi_k(y_0)R(y)$ for $y\ge y_0$ and $3\le k\le u<k+1$, where
\begin{dmath*}
\xi_k(y_0)\colonequals \left(\max_{t\ge y_0}\frac{1}{tR(t)}\right)y_0^{2-k}+C_7(y_0)\left(\frac{1}{(k-1)!}+\int_{k}^{k+1}\frac{dt}{\Gamma(t)}+\left(2\log 2-\frac{1}{2}\right)y_0^{2-k}+\int_{1}^{2}t^{-2}y_0^{t-k}\,dt+\int_{2}^{3}y_0^{t-k}\left(1-\log(t-1)+\frac{1}{(t-1)^2}\right)\,dt+2\int_{3}^{k}y_0^{t-k}\frac{dt}{\Gamma(t)}\right).
\end{dmath*}
As a direct corollary, we obtain
\begin{dmath*}
\sum_{k=2}^{\infty}\xi_{k}(y_0)=\frac{y_0}{y_0-1}\max_{t\ge y_0}\frac{1}{tR(t)}+C_7(y_0)\left(e-\frac{1}{2}+\int_{3}^{\infty}\frac{dt}{\Gamma(t)}+\frac{1}{y_0-1}\left(2\log 2-1+y_0I_{y_0}(2)+\int_{2}^{3}y_0^{t-2}\left(1-\log(t-1)+\frac{1}{(t-1)^2}\right)\,dt+2\int_{3}^{\infty}y_0^{\{t\}}\frac{dt}{\Gamma(t)}\right)\right),
\end{dmath*}
where we have applied partial summation to derive
\begin{align*}
\sum_{k=3}^{\infty}\int_{3}^{k}y_0^{t-k}\frac{dt}{\Gamma(t)}&=\left(\sum_{k=3}^{\infty}y_0^{-k}\right)\int_{3}^{\infty}y_0^{t}\frac{dt}{\Gamma(t)}-\int_{3}^{\infty}\left(\sum_{3\le k\le t}y_0^{-k}\right)y_0^{t}\frac{dt}{\Gamma(t)}\\
&=\frac{y_0^{-3}}{1-y_0^{-1}}\int_{3}^{\infty}y_0^{t}\frac{dt}{\Gamma(t)}-\int_{3}^{\infty}\frac{y_0^{t-3}(1-y_0^{-\lfloor t\rfloor+2})}{1-y_0^{-1}}\cdot\frac{dt}{\Gamma(t)}\\
&=\frac{1}{y_0-1}\int_{3}^{\infty}y_0^{\{t\}}\frac{dt}{\Gamma(t)}.
\end{align*}
For computational purposes, we can transform the last integral above by observing that
\[\int_{3}^{\infty}y_0^{\{t\}}\frac{dt}{\Gamma(t)}=\int_{0}^{1}\left(\sum_{n=0}^{\infty}\frac{1}{(t+2)\cdots(t+2+n)}\right)y_0^{t}\frac{dt}{\Gamma(t+2)}.\]
Let
\[\gamma(s,z)\colonequals\int_{0}^{z}v^{s-1}e^{-v}\,dv\]
be the lower incomplete gamma function, where $s\in\C$ with $\Re(s)>0$ and $z\ge0$. It is well known that
\[\gamma(s,z)=z^se^{-z}\sum_{n=0}^{\infty}\frac{z^n}{s(s+1)\cdots(s+n)},\]
from which it follows that \[\sum_{n=0}^{\infty}\frac{1}{(t+2)\cdots(t+2+n)}=\gamma(t+2,1)e.\]
Thus we obtain
\begin{dmath}\label{Equ:sumofxi_k}
\sum_{k=2}^{\infty}\xi_{k}(y_0)=\frac{y_0}{y_0-1}\max_{t\ge y_0}\frac{1}{tR(t)}+C_7(y_0)\left(e-\frac{1}{2}+\int_{3}^{\infty}\frac{dt}{\Gamma(t)}+\frac{1}{y_0-1}\left(2\log 2-1+y_0I_{y_0}(2)+\int_{2}^{3}y_0^{t-2}\left(1-\log(t-1)+\frac{1}{(t-1)^2}\right)\,dt+2e\int_{0}^{1}y_0^{t}\frac{\gamma(t+2,1)}{\Gamma(t+2)}\,dt\right)\right).
\end{dmath}
In Mathematica, the function $\gamma(t+2,1)$ can be evaluated by ``Gamma[t+2,0,1]".
\par Finally, we are ready to estimate $\eta(x,y)$. Let 
\[\eta_k(y)\colonequals \frac{1}{e^{\gamma}R(y)}\sup_{\substack{u\in[k,k+1)\\t\ge y}}|\eta(t^u,t)|\]
for $k\ge1$ and $y\ge y_0$, where the value of $\eta_1(y)$ is provided by (\ref{Equ:eta_1'}). Using (\ref{Equ:eta(x,y)1}) and the estimates for $E_3(h_k;y,u)$ with $y\ge y_0$ and $2\le k\le u<k+1$, we find
\[\eta_k(y)\le\eta_{k-1}(y)+\xi_k(y_0)\]
for all $k\ge2$ and $y\ge y_0$, from which we derive
\[\eta_k(y)\le\eta_1(y)+\sum_{\ell=2}^{k}\xi_{\ell}(y_0)\]
for all $k\ge1$ and $y\ge y_0$. Since $\eta_1(y)$ is decreasing on $[y_0,\infty)$, we have therefore shown that
\begin{equation}\label{Equ:eta(x,y)2}
|\eta(x,y)|\le e^{\gamma}\left(\eta_1(y_0)+\sum_{k=2}^{\infty}\xi_{k}(y_0)\right)R(y)
\end{equation}
for all $y\ge y_0$, where the infinite sum can be evaluated using (\ref{Equ:sumofxi_k}). To derive an explicit version of de Bruijn's result (\ref*{Equ:de Bruijn}), we observe that (\ref{Equ:Q(y)}), (\ref{Equ:eta(x,y)2}) and \cite[Theorem 23]{RS} imply that $Q(y)|\eta(x,y)|\le C_8(y_0)R(y)/\log y$ for all $y\ge y_0$, where
\[C_8(y_0)\colonequals\beta(y_0)\left(\eta_1(y_0)+\sum_{k=2}^{\infty}\xi_{k}(y_0)\right)\]
with
\[\beta(y_0)\colonequals\begin{cases}
~1,& \text{\hspace*{2mm}if~}3\le y_0<10^8,\\
~\exp(C_2(y_0)R(y_0)),&\text{\hspace*{2mm}if~}y_0\ge10^8.
\end{cases}\]
Hence, it follows that
\begin{equation}\label{Equ:ExplicitdeBruijn}
\left|\Phi(x,y)-\mu_y(u)e^{\gamma}x\log y\prod_{p\leq y}\left(1-\frac{1}{p}\right)\right|< \frac{C_8(y_0)xR(y)}{\log y}
\end{equation}
for all $y\ge y_0$. 

\section{Deduction of Theorem \ref{thm:explicit debruijn} and Corollary \ref{cor::explicit debruijn}}\label{S: Pf of Thm2&Cor3}
\par Now we apply (\ref{Equ:ExplicitdeBruijn}) to obtain explicit estimates for $\Phi(x,y)$ with specific choices of $R(y)$. Unconditionally, it has been shown \cite[Corollary 2]{MT} that
\[|\pi(z)-\li(z)|\le0.2593\frac{z}{(\log z)^{3/4}}\exp\left(-\sqrt{\frac{\log z}{6.315}}\right)\]
for all $z\ge229$. With $y_0\ge229$, the function
\[R(z)=0.2593(\log z)^{1/4}\exp\left(-\sqrt{\frac{\log z}{6.315}}\right)\]
is strictly decreasing on $[y_0,\infty)$ and satisfies (\ref*{Equ:pi(x)}) and (\ref*{Equ:intpi(x)}) with
\[C_0(y_0)=2\sqrt{\frac{6.315}{\log y_0}},\] 
since 
\begin{align*}
\int_{z}^{\infty}\frac{1}{t(\log t)^{3/4}}\exp\left(-\sqrt{\frac{\log t}{6.315}}\right)\,dt&=2\int_{\sqrt{\log z}}^{\infty}\frac{1}{\sqrt{t}}\exp\left(-\frac{t}{\sqrt{6.315}}\right)\,dt\\
&<\frac{2}{(\log z)^{1/4}}\int_{\sqrt{\log z}}^{\infty}\exp\left(-\frac{t}{\sqrt{6.315}}\right)\,dt\\
&=\frac{2\sqrt{6.315}}{(\log z)^{1/4}}\exp\left(-\sqrt{\frac{\log z}{6.315}}\right)
\end{align*}
for $z\ge y_0$. Numerical computation using Mathematica allows us to conclude that
\begin{equation}\label{Equ:ExplicitdeBruijn1}
\left|\Phi(x,y)-\mu_y(u)e^{\gamma}x\log y\prod_{p\leq y}\left(1-\frac{1}{p}\right)\right|< 4.403611\frac{x}{(\log y)^{3/4}}\exp\left(-\sqrt{\frac{\log y}{6.315}}\right)
\end{equation}
for all $x\ge y\ge 229$. Suppose now that $2\le y<229$. Using the inequalities $\Phi(x,y)<x/\log y$ \cite[Theorem]{F}, $\prod_{p\le y}(1-1/p)<e^{-\gamma}/\log y$ \cite[Theorem 23]{RS} and $0\le\mu_y(u)<1/\log y$, we have
\[\left|\Phi(x,y)-\mu_y(u)e^{\gamma}x\log y\prod_{p\leq y}\left(1-\frac{1}{p}\right)\right|<\frac{2x}{\log y}<4.403611\frac{x}{(\log y)^{3/4}}\exp\left(-\sqrt{\frac{\log y}{6.315}}\right)\]
for all $2\le y<229$. Combining this with (\ref{Equ:ExplicitdeBruijn1}) proves the first half of Theorem \ref{thm:explicit debruijn}. 
\par Under the assumption of the Riemann Hypothesis, 
it is known \cite[Corollary 1]{S} that
\[|\pi(z)-\li(z)|<\frac{1}{8\pi}\sqrt{z}\log z\]
for all $z\ge2{,}657$. With $y_0=2{,}657$ and
\[R(z)=\frac{\log^2z}{8\pi\sqrt{z}},\]
we have
\[\int_{z}^{\infty}\frac{|\pi(t)-\li(t)|}{t^2}\,dt\le \frac{1}{8\pi}\int_{z}^{\infty}\frac{\log t}{t^{3/2}}\,dt=\frac{\log z+2}{4\pi\sqrt{z}}\le C_0(y_0)R(z)\]
for $z\ge y_0$, where 
\[C_0(y_0)=\frac{2(\log y_0+2)}{\log^2 y_0}.\]
Therefore, we conclude by (\ref{Equ:ExplicitdeBruijn}) and numerical calculations that
\begin{equation}\label{Equ:ExplicitdeBruijn2}
\left|\Phi(x,y)-\mu_y(u)e^{\gamma}x\log y\prod_{p\leq y}\left(1-\frac{1}{p}\right)\right|< 0.184563\frac{x\log y}{\sqrt{y}}
\end{equation}
for all $x\ge y\ge 2{,}657$. The values of relevant constants are recorded in the table below.

\begin{table*}[ht]
\caption*{Table: Numerical Constants}
\begin{tabular}{| *{5}{c|} }
	\hline
	constants    & \multicolumn{2}{c|}{unconditional estimates}
	& \multicolumn{2}{c|}{conditional estimates}       \\ \hline
	
~$y_0$ &~229 &~$10^8$&~2{,}657&~$10^8$\\
~$R(y_0)$ &~.156576 &~ .097363&~.047992 &~.001351\\
~$C_0(y_0)$ &~2.156096 &~1.171019&~.317985 &~.120362\\
~$C_1(y_0)$ &~2.534430  &~1.279593&~.571800 &~.228936\\
~$C_2(y_0)$ &~2.548436  &~1.279593&~.575723 &~.228940\\
~$C_3(y_0)$ &~3.110976&~1.362717&~.579718 &~.228971\\
~$C_4(y_0)$ &~2.548436  &~1.279593&~.575723 &~.228940\\
~$C_5(y_0)$ &~3.131827&~1.362717&~.583750 &~.228975\\
~$C_6(y_0)$ &~2.534430 &~1.279593&~.571800 &~.228936\\
~$C_7(y_0)$ &~3.110976 &~1.362717&~.579718 &~.228971\\
~$C_8(y_0)$ &~16.982691 &~9.079975&~4.638553&~2.967998\\
~$\eta_1(y_0)$ &~6.236726&~3.628074&~2.697198 &~2.229726\\
~$\sum_{k=2}^{\infty}\xi_k(y_0)$ &~10.745960&~4.388310 &~1.941356 &~.737355\\\hline
\end{tabular}
\end{table*}
\FloatBarrier
\par To complete the proof of the second half of Theorem \ref{thm:explicit debruijn}, it remains to deal with the case $11\le y\le 2{,}657$. For simplicity of notation we set
\[D(x,y)\colonequals\Phi(x,y)-\mu_y(u)e^{\gamma}x\log y\prod_{p\leq y}\left(1-\frac{1}{p}\right).\]
Using Mathematica we find that
\begin{align*}
M&\colonequals\max\limits_{11\le z\le 2{,}657}\frac{\li(z)-\pi(z)}{\sqrt{z}\log z}<0.259141,\\
m&\colonequals\min\limits_{11\le z\le 2{,}657} e^{\gamma}\log z\prod_{p\le z}\left(1-\frac{1}{p}\right)>0.876248.
\end{align*}
If $\sqrt{x}\le y<x$, then 
\[\Phi(x,y)=\mu_y(u)x+(\pi(x)-\li(x))-(\pi(y)-\li(y))+1.\]
Note that $x\le y^2<10^{8}$. Since $\pi(z)<\li(z)$ for $2\le z\le10^{8}$ by \cite[Theorem 16]{RS} and  
\[\prod_{p\leq z}\left(1-\frac{1}{p}\right)<\frac{e^{-\gamma}}{\log z}\]
for $0<z\le10^{8}$ by \cite[Theorem 23]{RS}, we have
\begin{align}
|D(x,y)|&<(1-m)\left(1-y^{-1}\right)\frac{x}{\log y}+M\sqrt{x}\log x+1\nonumber\\
&\le\left((1-m)\left(1-y^{-1}\right)+M\frac{\log^2y}{\sqrt{y}}+\frac{\log y}{y}\right)\frac{x}{\log y},\label{Equ:1<=u<=2}
\end{align}
where we have used the fact that $\log x/\sqrt{x}$ is strictly decreasing on $[e^2,\infty)$. Numerical computation shows that the right side of (\ref{Equ:1<=u<=2}) is $<0.449774x\log y/\sqrt{y}$ for $11\le y\le 2{,}657$. Suppose now that $11\le y\le\sqrt{x}$. By \cite[Theorem 1]{FP}, Theorem \ref{thm:0.4LB} and \cite[Theorem 23]{RS} we have, for $11\le y\le2{,}657$,
\begin{align*}
D(x,y)&\le\left(0.6-\frac{m}{2}\left(1-y^{-1}\right)\right)\frac{x}{\log y}<0.449774\frac{x\log y}{\sqrt{y}},\\
D(x,y)&>(0.4-M_0)\frac{x}{\log y}>-0.449774\frac{x\log y}{\sqrt{y}}.
\end{align*}
This settles the case $11\le y\le 2{,}657$ and completes the proof of Theorem \ref{thm:explicit debruijn}.
\par The proof of Corollary \ref{cor::explicit debruijn} is similar, and we shall only sketch it. When $y\ge y_0$, where $y_0=229$ for the unconditional estimate and $y_0=2{,}657$ for the conditional estimate, we have by the triangle inequality that
\[|\Phi(x,y)-\mu_y(u)x|< |D(x,y)|+\left|1-e^{\gamma}\log y\prod_{p\leq y}\left(1-\frac{1}{p}\right)\right|\frac{x}{\log y}.\]
Then we bound $|D(x,y)|$ using the values of $C_8(y_0)$ listed in the table above. To estimate the second term, we use (\ref{Equ:Q(y)}) when $y\ge 10^8$ and the inequality 
\[m(y)<e^{\gamma}\log y\prod_{p\leq y}\left(1-\frac{1}{p}\right)<1\]
when $y_0\le y\le 10^8$, where $m(y)$ is given by
\[m(y)\colonequals\begin{cases}
	~0.983296,& \text{\hspace*{2mm}if~}229\le y\le 2{,}657,\\
	~0.996426,&\text{\hspace*{2mm}if~}2{,}657\le y<210{,}000,\\
	~0.999643,&\text{\hspace*{2mm}if~}210{,}000\le y\le10^8,
\end{cases}\]
according to \cite[Theorem 23]{RS} and Mathematica. This leads to the asserted bounds for $y\ge y_0$. Suppose now that $y\le y_0$. In this case, the proof of the unconditional bound is exactly the same as that of the unconditional bound in Theorem \ref{thm:explicit debruijn}. As for the conditional bound, we argue in the same way as in the proof of Theorem \ref{thm:explicit debruijn} to get
\[|\Phi(x,y)-\mu_y(u)x|\le\left(M\frac{\log^2y}{\sqrt{y}}+\frac{\log y}{y}\right)\frac{x}{\log y}\]
when $\sqrt{x}\le y<x$ and
\begin{align*}
|\Phi(x,y)-\mu_y(u)x|&\le\left(0.6-\frac{1}{2}\left(1-y^{-1}\right)\right)\frac{x}{\log y},\\
|\Phi(x,y)-\mu_y(u)x|&>(0.4-M_0)\frac{x}{\log y},
\end{align*}
when $11\le y\le\sqrt{x}$. Together, these inequalities yield the asserted conditional bound.

\begin{rmk}
The bounds in Theorem \ref{thm:explicit debruijn} and its corrollary may be improved. For example, the numerical values of the sum $\sum_{k=2}^{\infty}\xi_k(y_0)$ may be reduced by keeping $\rho$ (or even $|\omega'|$) in all of the relevant integrals, but of course the computational complexity is expected to increase as a cost. In addition, our method would allow an extension of the range $x\ge y\ge 11$ in the second half of Theorem \ref{thm:explicit debruijn} to the entire range $x\ge y\ge 2$ if we argue with $y_0=2{,}657$ replaced by some smaller value and enlarge the constant 0.449774.
\end{rmk}

{\noindent\bf Acknowledgment.} The author would like to thank his advisor C. Pomerance for his helpful comments and suggestions.

\medskip

\bibliographystyle{amsplain}
\bibliography{bibliography}

\begin{bibdiv}
\begin{biblist}

\bib{BD}{article}{
      author={Berkane, D.},
      author={Dusart, P.},
       title={{On a constant related to the prime counting function}},
        date={2016},
     journal={Mediterr. J. Math.},
      volume={13},
      number={3},
       pages={929\ndash 938},
      review={\MR{3513147}},
}

\bib{B}{article}{
      author={Broadbent, S.},
      author={Kadiri, H.},
      author={Lumley, A.},
      author={Ng, N.},
      author={Wilk, K.},
       title={{Sharper bounds for the Chebyshev function $\theta(x)$}},
        date={2021},
     journal={Math. Comp.},
      volume={90},
      number={331},
       pages={2281\ndash 2315},
      review={\MR{4280302}},
}

\bib{Bu}{article}{
      author={Buchstab, A.~A.},
       title={{Asymptotic estimates of a general number-theoretic function
  (Russian, German summary)}},
        date={1937},
     journal={Mat. Sb.},
      volume={44},
       pages={1239\ndash 1246},
}

\bib{Br}{article}{
      author={de~Bruijn, N.~G.},
       title={{On the number of uncancelled elements in the sieve of
  Eratosthenes}},
        date={1950},
     journal={Nederl. Akad. Wetensch., Proc.},
      volume={53},
       pages={803\ndash 812},
      review={\MR{0035785}},
}

\bib{D}{article}{
      author={Dusart, P.},
       title={{Explicit estimates for some functions over primes}},
        date={2018},
     journal={Ramanujan J.},
      volume={45},
      number={1},
       pages={227\ndash 251},
      review={\MR{3745073}},
}

\bib{F}{article}{
      author={Fan, K.~(S.)},
       title={{An inequality for the distribution of numbers free of small
  prime factors}},
        date={2022},
     journal={Integers},
      volume={22},
      number={A26},
       pages={12 pp.},
      review={\MR{4390731}},
}

\bib{FP}{article}{
      author={Fan, K.~(S.)},
      author={Pomerance, C.},
       title={{An inequality related to the sieve of Eratosthenes}},
        date={2023},
     journal={J. Number Theory},
       pages={to appear},
}

\bib{JR}{article}{
      author={Jurkat, W.~B.},
      author={Richert, H.-E.},
       title={{An improvement of Selberg’s sieve method. I}},
        date={1965},
     journal={Acta Arith.},
      volume={11},
       pages={217\ndash 410},
      review={\MR{0202680}},
}

\bib{MT}{article}{
      author={Mossinghoff, M.~J.},
      author={Trudgian, T.~S.},
       title={{Nonnegative trigonometric polynomials and a zero-free region for
  the Riemann zeta-function}},
        date={2015},
     journal={J. Number Theory},
      volume={157},
       pages={329\ndash 349},
      review={\MR{3373245}},
}

\bib{Rama}{article}{
      author={Rama\'re, O.},
       title={{An explicit density estimate for Dirichlet $L$-series}},
        date={2016},
     journal={Math. Comp.},
      volume={85},
      number={297},
       pages={325\ndash 356},
      review={\MR{3404452}},
}

\bib{RS}{article}{
      author={Rosser, J.~B.},
      author={Schoenfeld, L.},
       title={{Approximate formulas for some functions of prime numbers}},
        date={1962},
     journal={Illinois J. Math.},
      volume={6},
       pages={64\ndash 94},
      review={\MR{0137689}},
}

\bib{S}{article}{
      author={Schoenfeld, L.},
       title={{Sharper bounds for the Chebyshev functions $\theta(x)$ and
  $\psi(x)$. II}},
        date={1976},
     journal={Math. Comp.},
      volume={30},
      number={134},
       pages={337\ndash 360},
      review={\MR{0457374}},
}

\bib{T}{book}{
      author={Tenenbaum, G.},
       title={Introduction to analytic and probabilistic number theory},
      series={3rd. ed., Grad. Stud. Math.},
   publisher={Amer. Math. Soc.},
     address={Providence, RI.},
        date={2015},
      volume={163},
      review={\MR{3363366}},
}

\end{biblist}
\end{bibdiv}
\end{document}